\documentstyle{article}

\def\G{$\bf{\Gamma}$}
\def\cirk{\,{\raisebox{.3ex}{\tiny $\circ$}}\,}
\def\kon{\wedge}
\def\eL{\mbox{$\cal L$}}

\def\koc{\begin{picture}(10,7)\unitlength1pt
\put(1.5,0){\line(0,1){7}} \put(8.5,0){\line(0,1){7}}
\put(1.5,0){\line(1,0){7}} \put(1.5,7){\line(1,0){7}}
\end{picture}}

\def\koci{\begin{picture}(7,5)
\put(1,0){\line(0,1){5}} \put(6,0){\line(0,1){5}}
\put(1,0){\line(1,0){5}} \put(1,5){\line(1,0){5}}
\end{picture}}

\def\eI{\mbox{$\cal I$}}
\def\prop#1#2{\vspace{2ex} \noindent{\sc #1.} {\it #2} \par \vspace{2ex}}
\def\dkz{\noindent{\sc Proof. }}
\def\qed{\hfill $\dashv$}
\def\df{\mbox{\scriptsize{\it df}}}

\def\mak{{\bf\^ A}}
\def\g22{\gamma^\rightarrow_{{\bf 2},{\bf 2}}}

\begin{document}

\title{Associativity as Commutativity}
\author{{\sc Kosta Do\v sen} and {\sc Zoran Petri\' c}
\\[1ex]
{\small Mathematical Institute, SANU}\\[-.5ex]
{\small Knez Mihailova 35, p.f. 367, 11001 Belgrade,
Serbia}\\[-.5ex]
{\small email: \{kosta, zpetric\}@mi.sanu.ac.yu}}
\date{}
\maketitle

\vspace{-6ex}

\begin{abstract}
\noindent It is shown that coherence conditions for monoidal
categories concerning associativity are analogous to coherence
conditions for symmetric or braided strictly monoidal categories,
where associativity arrows are identities. Mac Lane's pentagonal
coherence condition for associativity is decomposed into
conditions concerning commutativity, among which we have a
condition analogous to naturality and a degenerate case of Mac
Lane's hexagonal condition for commutativity. This decomposition
is analogous to the derivation of the Yang-Baxter equation from
Mac Lane's hexagon and the naturality of commutativity. The
pentagon is reduced to an inductive definition of a kind of
commutativity.
\end{abstract}

\vspace{.3cm}

\noindent {\small {\it Mathematics Subject Classification} ({\it
2000}): 18D10, 18A05}

\vspace{.5ex}

\noindent {\small {\it Keywords$\,$}: monoidal categories,
symmetric monoidal categories, braided monoidal categories,
coherence, Mac Lane's pentagon, Mac Lane's hexagon, insertion}

\vspace{.5ex}

\noindent {\small {\it Acknowledgement$\,$}: We would like to
thank Slobodan Vujo\v sevi\' c and an anonymous referee for
reading the paper and making helpful suggestions. We are grateful
to Milo\v s Ad\v zi\' c for drawing our attention to \cite{M06}.}

\vspace{1cm}

\vspace{-7ex}

\baselineskip=1.25\baselineskip

\section{Introduction}
Associativity is a kind of commutativity. To see why, conceive of
${(a\cdot(c\cdot b))}$ as
${(a\cdot\underline{\;\;\;}\,)\circ(\,\underline{\;\;\;}\cdot b)}$
applied to $c$. We have

\vspace{-1ex}

\begin{tabbing}
\hspace{7em}$((a\cdot\underline{\;\;\;}\,)\circ(\,\underline{\;\;\;}\cdot
b))(c)\;$\=$=(a\cdot\underline{\;\;\;}\,)((\,\underline{\;\;\;}\cdot
b)(c))$
\\*[.5ex]
\>$=(a\cdot\underline{\;\;\;}\,)((c\cdot b))$
\\*[.5ex]
\>$=(a\cdot(c\cdot b)).$

\end{tabbing}

\vspace{-1ex}

\noindent Then associativity goes from
${(a\cdot\underline{\;\;\;}\,)\circ(\,\underline{\;\;\;}\cdot b)}$
to ${(\,\underline{\;\;\;}\cdot
b)\circ(a\cdot\underline{\;\;\;}\,)}$, which when applied to $c$
yields ${((a\cdot c)\cdot b)}$.

The purpose of this paper is to exploit this idea to show that
monoidal categories may be conceived as a kind of symmetric
strictly monoidal categories, where associativity arrows are
identities. As a matter of fact, this analogy holds also with
braided strictly monoidal categories. More precisely, we show that
coherence conditions for monoidal categories concerning
associativity are analogous to coherence conditions concerning
commutativity (i.e.\ symmetry or braiding) for symmetric (see
\cite{ML63}, \cite{ML71} or \cite{DP04}) or braided (see
\cite{JS93} and \cite{ML71}, second edition, Chapter IX) strictly
monoidal categories. In particular, Mac Lane's pentagonal
coherence condition for associativity (see Section 4 below) is
decomposed into conditions concerning commutativity, among which
we have a condition analogous to naturality and degenerate cases
of Mac Lane's hexagonal condition for commutativity. (The hexagon
becomes a triangle, because associativity arrows are identities,
or a two-sided figure.) This decomposition is analogous to the
derivation of the Yang-Baxter equation from Mac Lane's hexagon and
the naturality of commutativity (see Section 5 below).

To achieve that, we replace the algebra freely generated with one
binary operation (denoted by $\kon$) by an isomorphic algebra
generated with a family of partial operations we call
\emph{insertion} (denoted by $\triangleleft_n$); insertion is
analogous to the composition $\circ$ at the beginning of this
text, or to functional application. (This procedure is like
Achilles' introduction of $\alpha$ in \cite{D98}.) The latter
algebra is more complicated, and is not free any more, but it
enables us to present associativity arrows as commutativity
arrows.

In the next section we state precisely these matters concerning
insertion. After that we introduce a category \G, which in the
remainder of the paper is shown isomorphic to a free monoidal
category without unit. In \G\ the associativity arrows appear as a
kind of commutativity arrows, and coherence conditions for \G\
take the form of an inductive definition of these commutativity
arrows. Our decomposition of Mac Lane's pentagon is in the last
section.

We work with categories where associativity is an isomorphism,
because this is the standard approach, but our treatment is easily
transferred to categories where associativity arrows are not
necessarily isomorphisms, which in \cite{DP04} (Section 4.2) are
called \emph{semiassociative} categories (see also \cite{L72}). As
monoidal categories, semiassociative categories are coherent in
Mac Lane's ``all diagrams commute'' sense. With semiassociative
categories, the commutativity corresponding to associativity only
ceases to be an isomorphism, and all the rest remains as in the
text that follows.

Among the coherence conditions for the main kinds of categories
with structure, Mac Lane's pentagon seems to be more mysterious
than the others. Our decomposition of the pentagon goes towards
dispelling the mystery. The pentagon is reduced to an inductive
definition of a kind of commutativity.

If the associativity arrows are isomorphisms, then the pentagon
yields the definition of an associativity arrow complex in one of
its indices in terms of associativity arrows simpler in that
index, but more complex in the other two indices. There is no
reduction of complexity in all the indices, and no real inductive
definition. Our approach, which works also in the absence of
isomorphism, as we said above, gives a real inductive definition.

\section{Insertion}
Let $\eL_1$ be the set of words (finite sequences of symbols) in
the alphabet ${\{\koc,\kon,(,)\}}$ defined inductively by

\begin{tabbing}

\hspace{11em}\=$\koc\in\eL_1$,
\\*[1ex]
\> if $X,Y\in\eL_1$, then $(X\kon Y)\in\eL_1$.

\end{tabbing}

Let $\eL_2$ be ${\eL_1-\{\koc\}}$. The elements of $\eL_2$ may be
identified with finite planar binary trees with more than one
node, while $\koc$ is the trivial one-node tree. In this section,
we use $X$, $Y$ and $Z$ for the members of $\eL_1$. (Starting from
the end of the section, we change this notation to $A$, $B$, $C,
\ldots$) We omit the outermost pair of parentheses of the members
of $\eL_1$, taking them for granted. We make the same omission in
other analogous situations later on.

Let ${\bf N}^+$ be the set of natural numbers greater than $0$,
and let $\eL'$ be the set of words in the alphabet ${\{{\bf
1},{\bf 2}\}\cup\{\triangleleft_n\,|\,n\in{\bf N}^+\}}$ defined
inductively by the following clauses that involve also an
inductive definition of a map $|\;|$ from $\eL'$ to ${\bf N}^+$:

\begin{tabbing}

\hspace{13.5em}\= ${\bf 1}\in\eL'$ and $|{\bf 1}|=1$,
\\*[1ex]
\> ${\bf 2}\in\eL'$ and $|{\bf 2}|=2$,
\\[1.5ex]
\hspace{1em} if $A,B\in\eL'$ and $1\leq n\leq |A|$, then \\[.5ex]
\> $(A\triangleleft_n B)\in\eL'$ and $|(A\triangleleft_n
B)|=|A|+|B|-1$.

\end{tabbing}

Let $\eL''$ be defined as $\eL'$ save that we omit the first
clause above involving $\bf 1$, and we replace $\eL'$ by $\eL''$
in the two remaining clauses. For the members of $\eL'$ we use
$A$, $B$, $C,\ldots$, sometimes with indices. As we did for
$\eL_1$, we omit the outermost pair of parentheses of the members
of $\eL'$.

We define the equational calculus $\eI''$ in $\eL''$ (i.e.\ a
calculus whose theorems are equations between members of $\eL''$)
by assuming reflexivity, symmetry and transitivity of equality,
the rule that if ${A=B}$ and ${C=D}$, then ${A\triangleleft_n
C=B\triangleleft_n D}$, provided ${A\triangleleft_n C}$ and
${B\triangleleft_n D}$ are defined, and the two axioms

\begin{tabbing}

\hspace{1em}\=(\emph{assoc}~1) \quad \= $(A\triangleleft_n
B)\triangleleft_m C = A\triangleleft_n (B\triangleleft_{m-n+1} C)$
\quad\quad\= if $n\leq m<n+|B|$,
\\[1.5ex]
\>(\emph{assoc}~2)\>$(A\triangleleft_n B)\triangleleft_m C
=(A\triangleleft_{m-|B|+1} C)\triangleleft_n B$ \> if $n+|B|\leq
m$.

\end{tabbing}

Note that the condition ${n\leq m<n+|B|}$ in
\mbox{(\emph{assoc}~1)} follows from the legitimacy of
${B\triangleleft_{m-n+1} C}$. Note also that in both
\mbox{(\emph{assoc}~1)} and \mbox{(\emph{assoc}~2)} we have
${n\leq m}$. The equation \mbox{(\emph{assoc}~2)} could be
replaced by

\begin{tabbing}

\hspace{1em}\=(\emph{assoc}~1) \quad \= $(A\triangleleft_n
B)\triangleleft_m C = A\triangleleft_n (B\triangleleft_{m-n+1} C)$
\quad\quad\= if $n\leq m<n+|B|$\kill

\>\>$(A\triangleleft_n B)\triangleleft_m C = (A\triangleleft_m
C)\triangleleft_{n+|C|-1} B$\> if $m<n$.

\end{tabbing}
(The equations (\emph{assoc}~1) and (\emph{assoc}~2) are analogous
to the two associativity equations for the cut operation one finds
in multicategories; see \cite{L69} and \cite{L89}, Section~3.
Analogous equations are also found in the definition of operad;
see \cite{M06}, Section~1.)

The equational calculus $\eI'$ in $\eL'$ is defined as $\eI''$
with the additional axiom

\begin{tabbing}

\hspace{1em}\=(\emph{assoc}~1) \quad \= $(A\triangleleft_n
B)\triangleleft_m C = A\triangleleft_n (B\triangleleft_{m-n+1} C)$
\quad\quad\= if $n\leq m<n+|B|$\kill

\>(\emph{unit})\>${\bf 1}\triangleleft_1 A=A\triangleleft_n {\bf
1}=A$

\end{tabbing}

\noindent(whose analogue one also finds in multicategories). Our
purpose now is to interpret $\eL'$ in $\eL_1$. This will make
clear the meaning of the axioms of $\eI'$.

For $X$ in $\eL_1$, let ${|X|}$ be the number of occurrences of
$\koc$ in $X$. We define in $\eL_1$ the partial operation of
\emph{insertion} $\unlhd_n$ by the following inductive clauses:

\begin{tabbing}

\hspace{4em}\=$\koc\unlhd_1 Z=Z$,
\\*[2ex]
\>$(X\kon Y)\unlhd_n Z=\left\{\begin{array}{ll} (X\unlhd_n Z)\kon
Y & \mbox{\rm{if }} 1\leq n\leq |X|
\\[1.5ex]
X\kon(Y\unlhd_{n-|X|} Z) & \mbox{\rm{if }} |X|<n\leq |X|+|Y|.
\end{array}\right.$

\end{tabbing}

We define insertion in $\eL_2$ by replacing the clause
$\koc\unlhd_1 Z=Z$ above by the clauses

\[
\begin{array}{l}

(\koc\kon\koc)\unlhd_1 Z=Z\kon\koc,
\\[1.5ex]
(\koc\kon\koc)\unlhd_2 Z=\koc\kon Z.

\end{array}
\]

Insertion gets its name from the fact that ${X\unlhd_n Z}$ is
obtained by \emph{inserting} $Z$ at the place of the $n$-th
occurrence of $\koc$ in $X$, starting from the left; namely, the
$n$-th leaf of the tree corresponding to $X$ becomes the root of
the tree corresponding to $Z$, and the resulting tree corresponds
to ${X\unlhd_n Z}$. Insertion is called \emph{grafting} in
\cite{S97}, and particular instances of insertion, which one finds
in the source and target of the arrows $\gamma^\rightarrow_{A,B}$
in Section 3 below, are called \emph{under} and \emph{over} in
\cite{L02} (Section 1.5).

We interpret $\eL''$ in $\eL_2$, i.e., we define a function $v$
from $\eL''$ to $\eL_2$, in the following manner:

\[
\begin{array}{l}

v({\bf 2})=\koc\kon\koc,
\\[1.5ex]
v(A\triangleleft_n B)=v(A)\unlhd_n v(B).

\end{array}
\]

\noindent For this definition to be correct, we must check that
${|A|=|v(A)|}$, which is easily done by induction on the length of
${|A|}$.

We prove first the following by an easy induction on the length of
derivation.

\prop{Soundness}{If ${A=B}$ in $\eI''$, then ${v(A)=v(B)}$.}

\noindent Our purpose is to prove also the converse:

\prop{Completeness}{If ${v(A)=v(B)}$, then ${A=B}$ in $\eI''$.}

For every $A$ in $\eL''$ we define the natural number ${c(A)}$
inductively as follows:

\[
\begin{array}{l}

c({\bf 2})=2,
\\[1.5ex]
c(B\triangleleft_n C)=c(B)(c(C)+1).

\end{array}
\]

\noindent Let ${s(A)}$ be the sum of the indices $n$ of all the
occurrences of $\triangleleft_n$ in $A$, and let
${d(A)=c(A)+s(A)}$. Then we can easily check that if ${A=B}$ is an
instance of \mbox{(\emph{assoc}~1)} or \mbox{(\emph{assoc}~2)},
then ${d(A)>d(B)}$.

Let a member of $\eL''$ be called \emph{normal} when it has no
part of the form of the left-hand side of \mbox{(\emph{assoc}~1)}
or \mbox{(\emph{assoc}~2)}, i.e.\ no part of the form
${(A\triangleleft_n B)\triangleleft_m C}$ for ${n\leq m}$. It can
be shown that a normal member of $\eL''$ is of one of the
following forms:

\[
({\bf 2}\triangleleft_2 A_2)\triangleleft_1 A_1, \;\;\;{\bf
2}\triangleleft_2 A_2, \;\;\;{\bf 2}\triangleleft_1 A_1,
\;\;\;{\bf 2},
\]

\noindent for $A_1$ and $A_2$ normal. These are the four
\emph{normal types}.

Then it is easy to show by applying \mbox{(\emph{assoc}~1)} and
\mbox{(\emph{assoc}~2)} from left to right that for every $A$ in
$\eL''$ there is a normal $A'$ such that ${A=A'}$ in $\eI''$. We
can also show the following.

\prop{Auxiliary Lemma}{If $A$ and $B$ are normal and
${v(A)=v(B)}$, then $A$ and $B$ coincide.}

\dkz If ${v(A)=v(B)}$, then $A$ and $B$ must be of the same normal
type (otherwise, clearly, ${v(A)\neq v(B)}$). If $A$ is ${({\bf
2}\triangleleft_2 A_2)\triangleleft_1 A_1}$ and $B$ is ${({\bf
2}\triangleleft_2 B_2)\triangleleft_1 B_1}$, then we conclude that
${v(A_1)=v(B_1)}$ and ${v(A_2)=v(B_2)}$, and we reason by
induction. We reason analogously for the second and third normal
type, and the normal type ${\bf 2}$ provides the basis of the
induction.\qed

\vspace{2ex}

To prove Completeness, suppose ${v(A)=v(B)}$. Let ${A=A'}$ and
${B=B'}$ in $\eI''$ for $A'$ and $B'$ normal. Then by Soundness we
have ${v(A')=v(A)=v(B)=v(B')}$, and so, by the Auxiliary Lemma,
$A'$ and $B'$ coincide. It follows that ${A=B}$ in $\eI''$, which
proves Completeness.

We can now also show that if ${A=A'}$ and ${A=A''}$ in $\eI''$ for
${A'}$ and ${A''}$ normal, then ${A'}$ and ${A''}$ coincide. This
follows from Soundness and the Auxiliary Lemma. (One could show
this uniqueness of normal form directly in $\eL''$, without
proceeding via $v$ and $\eL_2$, by relying on confluence
techniques, as in the lambda calculus or term-rewriting systems.
In such a proof, diagrams analogous to Mac Lane's pentagon and the
Yang-Baxter equation would arise.)

We interpret $\eL'$ in $\eL_1$ by extending the definition of $v$
from $\eL''$ to $\eL_2$ with the clause ${v({\bf 1})=\koc}$. Then
we can prove Soundness and Completeness with $\eI''$ replaced by
$\eI'$. In reducing a member of $\eL'$ to a normal member of
$\eL''$ or to ${\bf 1}$ we get rid first of all superfluous
occurrences of ${\bf 1}$, by relying on the equations
\mbox{(\emph{unit})}. Otherwise, the proof proceeds as before.

We can factorize $\eL'$ through the smallest equivalence relation
such that the equations of $\eI'$ are satisfied, and obtain a set
of equivalence classes isomorphic to $\eL_1$. For the equivalence
classes $[A]$ and $[B]$ we define $\kon$ by

\[
[A]\kon[B]=_{\df}\; [({\bf 2}\triangleleft_2 B)\triangleleft_1 A],
\]

\noindent and the isomorphism $i$ from $\eL_1$ to $\eL'$ is
defined by

\[
\begin{array}{l}
i(\koc)=[{\bf 1}],
\\[1.5ex]
i(X\kon Y)=i(X)\kon i(Y).

\end{array}
\]

\noindent The inverse $i^{-1}$ of $i$ is defined by
${i^{-1}([A])=v(A)}$. (To verify that $i$ and $i^{-1}$ are inverse
to each other, we rely on the fact that every ${[A]}$ is equal to
${[A']}$ for $A'$ being normal or ${\bf 1}$.)

We designate the equivalence class $[A]$ by $A$, and to designate
the elements of $\eL_1$ we can then use the notation introduced
for $\eL'$. This means that we can write $A$, $B$, $C,\ldots$ for
$X$, $Y$, $Z,\ldots$, we can write ${\bf 2}$ for ${\koc\kon\koc}$,
and we can write $\triangleleft_n$ for $\unlhd_n$. We have for
$\eL_1$ the equation

\[
A\kon B=({\bf 2}\triangleleft_2 B)\triangleleft_1 A.
\]

\noindent We will write however $\koc$ instead of {\bf 1}, and
reserve {\bf 1} with a subscript for the name of an arrow.

\section{The category \G}
The objects of the category \G\ are the elements of $\eL_1$. To
define the arrows of \G, we define first inductively the
\emph{arrow terms} of \G\ in the following way:

\[
\begin{array}{l}
{\bf 1}_A\!:A\rightarrow A,
\\[1ex]
\gamma^\rightarrow_{A,B}\!:A\triangleleft_{|A|} B\rightarrow
B\triangleleft_1 A,
\\[1.5ex]
\gamma^\leftarrow_{A,B}\!:B\triangleleft_1 A\rightarrow
A\triangleleft_{|A|} B
\end{array}
\]

\noindent are arrow terms of \G\ for all objects $A$ and $B$; if
${f\!:A\rightarrow B}$ and ${g\!:C\rightarrow D}$ are arrow terms
of \G, then ${g\cirk f\!:A\rightarrow D}$ is an arrow term of \G,
provided $B$ is $C$, and ${f\triangleleft_n g\!: A\triangleleft_n
C\rightarrow B\triangleleft_n D}$ is an arrow term of \G, provided
${1\leq n\leq |A|}$ and ${1\leq n\leq |B|}$. Note that for all
arrow terms ${f\!:A\rightarrow B}$ of \G\ we have ${|A|=|B|}$; we
write $|f|$ for $|A|$, which is equal to $|B|$.

The \emph{arrows} of \G\ are equivalence classes of arrow terms of
\G\ (cf.\ \cite{DP04}, Section 2.3) such that the following
equations are satisfied:

\begin{tabbing}

\hspace{1em}\=(\emph{cat}~1)\hspace{4em}\=${\bf 1}_B\cirk
f=f\cirk{\bf 1}_A=f$, \quad\quad for $f\!:A\rightarrow B$,
\\*[1.5ex]
\>(\emph{cat}~2)\>$(h\cirk g)\cirk f=h\cirk(g\cirk f)$,
\\[2ex]
\>(\emph{bif}~1)\>${\bf 1}_A\triangleleft_n {\bf 1}_B ={\bf
1}_{A\triangleleft_n B}$,
\\*[1.5ex]
\>(\emph{bif}~2)\>$(f_2\cirk f_1)\triangleleft_n(g_2\cirk
g_1)=(f_2\triangleleft_n g_2)\cirk(f_1\triangleleft_n g_1)$,
\\[2ex]
for $1\leq n\leq |f|$ and $n\leq m\leq |f|+|g|-1$
\\*[1.5ex]
\>(\emph{assoc}~1$\rightarrow$)\>$(f\triangleleft_n
g)\triangleleft_m h = f\triangleleft_n (g\triangleleft_{m-n+1} h)$
\quad\quad\= if $n\leq m<n+|g|$,
\\*[1.5ex]
\>(\emph{assoc}~2$\rightarrow$)\>$(f\triangleleft_n
g)\triangleleft_m h =(f\triangleleft_{m-|g|+1} h)\triangleleft_n
g$ \> if $n+|g|\leq m$,
\\[1.5ex]
\>(\emph{unit}~$\rightarrow$)\>${\bf 1}_{\koci}\triangleleft_1
f=f\triangleleft_n {\bf 1}_{\koci}=f$,
\\[2ex]
\>($\gamma$~\emph{nat})\>$\gamma^\rightarrow_{B,D}\cirk(f\triangleleft_{|A|}
g)=(g\triangleleft_1 f)\cirk\gamma^\rightarrow_{A,C}$,
\\[1.5ex]
\>($\gamma\gamma$)\>$\gamma^\leftarrow_{A,B}\cirk\gamma^\rightarrow_{A,B}={\bf
1}_{A\triangleleft_{|A|} B}$,\quad\quad
$\gamma^\rightarrow_{A,B}\cirk\gamma^\leftarrow_{A,B}={\bf
1}_{B\triangleleft_1 A}$,
\\[2ex]
\>($\gamma${\bf
1})\>$\gamma^\rightarrow_{\koci,A}=\gamma^\rightarrow_{A,\koci}=
{\bf 1}_A$,
\\[2ex]
\>(\emph{hex}~1)\>$\gamma^\rightarrow_{A\triangleleft_{|A|}
B,C}=(\gamma^\rightarrow_{A,C}\triangleleft_{|A|}{\bf
1}_B)\cirk({\bf 1}_A\triangleleft_{|A|}\gamma^\rightarrow_{B,C})$,
\\*[1.5ex]
\>(\emph{hex}~1\emph{a})\>$\gamma^\rightarrow_{A\triangleleft_n
B,C}=\gamma^\rightarrow_{A,C}\triangleleft_n{\bf
1}_B$\hspace{5em}\= if $1\leq n<|A|$,
\\[2ex]
\>(\emph{hex}~2)\>$\gamma^\rightarrow_{C,A\triangleleft_1 B}=({\bf
1}_A\triangleleft_1\gamma^\rightarrow_{C,B})\cirk
(\gamma^\rightarrow_{C,A}\triangleleft_{|C|}{\bf 1}_B)$,
\\*[1.5ex]
\>(\emph{hex}~2\emph{a})\>$\gamma^\rightarrow_{C,A\triangleleft_n
B}=\gamma^\rightarrow_{C,A}\triangleleft_{n+|C|-1}{\bf 1}_B$\> if
$1<n\leq |A|$.

\end{tabbing}

\noindent We also assume besides reflexivity, symmetry and
transitivity of equality that if ${f=g}$ and ${h=j}$, then for
$\alpha$ being $\cirk$ or $\triangleleft_n$ we have ${f\alpha
h=g\alpha j}$, provided ${f\alpha h}$ and ${g\alpha j}$ are
defined.

The equations \mbox{(\emph{cat}~1)} and \mbox{(\emph{cat}~2)} make
of \G\ a category. The equations \mbox{(\emph{bif}~1)} and
\mbox{(\emph{bif}~2)} are analogous to bifunctorial equations. The
equations \mbox{(\emph{assoc}~1$\rightarrow$)},
\mbox{(\emph{assoc}~2$\rightarrow$)} and
\mbox{(\emph{unit}~$\rightarrow$)} are analogous to naturality
equations. In \mbox{(\emph{assoc}~1$\rightarrow$)} the
associativity arrows with respect to $\triangleleft_n$ are not
written down because they are identity arrows, in virtue of the
equation \mbox{(\emph{assoc}~1)} on objects. Analogous remarks
hold for \mbox{(\emph{assoc}~2$\rightarrow$)} and
\mbox{(\emph{unit}~$\rightarrow$)}. The equation
\mbox{($\gamma$~\emph{nat})} is analogous to a naturality
equation, and the equations \mbox{($\gamma\gamma$)} say that
$\gamma^\rightarrow_{A,B}$ is an isomorphism, with inverse
$\gamma^\leftarrow_{A,B}$.

The equation \mbox{($\gamma{\bf 1}$)} is auxiliary, and would not
be needed if we had assumed $\gamma^\rightarrow_{A,B}$ and
$\gamma^\leftarrow_{A,B}$ only for $A$ and $B$ different from
$\koc$. The equations \mbox{(\emph{hex}~1)} and
\mbox{(\emph{hex}~2)} are analogous to Mac Lane's hexagonal
equation of symmetric monoidal categories (see \cite{ML63},
\cite{ML71}, Section VII.7, or \cite{DP04}, Section~5.1). Here the
associativity arrows with respect to $\triangleleft_n$ are
identity arrows, in virtue of the equation \mbox{(\emph{assoc}~1)}
on objects (and so instead of hexagons we have triangles; cf.\ the
equation \mbox{($c$~\emph{hex}~1)} in Section 5 below). Finally,
the equations \mbox{(\emph{hex}~1\emph{a})} and
\mbox{(\emph{hex}~2\emph{a})}, together with \mbox{($\gamma{\bf
1}$)}, \mbox{(\emph{hex}~1)} and \mbox{(\emph{hex}~2)}, enable us
to define inductively $\gamma^\rightarrow_{A,B}$ for all $A$ and
$B$ in terms of the identity arrows ${\bf 1}_A$, the arrows

\vspace{-.5ex}

\[
\g22\!: \koc\kon(\koc\kon\koc)\rightarrow(\koc\kon\koc)\kon\koc
\]

\noindent and the operations on arrows $\cirk$ and
$\triangleleft_n$. Relying on ($\gamma\gamma$), we can proceed
analogously for $\gamma^\leftarrow_{A,B}$ by using instead of
$\g22$ the arrows

\[
\gamma^\leftarrow_{{\bf 2},{\bf 2}}\!:
(\koc\kon\koc)\kon\koc\rightarrow \koc\kon(\koc\kon\koc).
\]

\noindent The equations \mbox{(\emph{hex}~1\emph{a})} and
\mbox{(\emph{hex}~2\emph{a})} are also analogous to Mac Lane's
hexagon mentioned above (due to the presence of (\emph{assoc}~2)
too, the collapse is however not any more into a triangle, but
into a two-sided figure).

\section{The category \mak}

The category \mak\ has the same objects as \G; namely, the
elements of $\eL_1$. To define the arrows of \mak, we define first
inductively the \emph{arrow terms} of \mak\ in the following way:

\vspace{-.5ex}

\[
\begin{array}{l}

{\bf 1}_A\!:A\rightarrow A,
\\[1ex]
b^\rightarrow_{A,B,C}\!:A\kon(B\kon C)\rightarrow(A\kon B)\kon C,
\\*[1.5ex]
b^\leftarrow_{A,B,C}\!:(A\kon B)\kon C\rightarrow A\kon(B\kon C)
\end{array}
\]

\noindent are arrow terms of \mak\ for all objects $A$, $B$ and
$C$; if ${f\!:\!A\rightarrow B}$ and ${g\!:\!C\rightarrow D}$ are
arrow terms of \mak, then ${g\cirk f\!:\!A\rightarrow D}$ is an
arrow term of \mak, provided $B$ is $C$, and ${f\kon g\!:\!A\kon
C\rightarrow B\kon D}$ is an arrow term of \mak.

The \emph{arrows} of \mak\ are equivalence classes of arrow terms
of \mak\ such that the following equations are satisfied:
\mbox{(\emph{cat}~1)}, \mbox{(\emph{cat}~2)},
\mbox{(\emph{bif}~1)} and \mbox{(\emph{bif}~2)} with
$\triangleleft_n$ replaced by $\kon$, and moreover

\begin{tabbing}
\mbox{\hspace{2em}}\= \mbox{($b$~\emph{nat})} \hspace{1em}\=
$b^\rightarrow_{B,D,F}\cirk(f\kon(g\kon h))=((f\kon g)\kon h)\cirk
b^\rightarrow_{A,C,E}$,
\\[2ex]
\>  \mbox{($bb$)}\> $b^\leftarrow_{A,B,C}\cirk
b^\rightarrow_{A,B,C}={\bf 1}_{A\kon(B\kon C)}$,\quad\quad
$b^\rightarrow_{A,B,C}\cirk b^\leftarrow_{A,B,C}={\bf 1}_{(A\kon
B)\kon C}$,
\\[2ex]
\> \mbox{($b$5)} \> $b^\rightarrow_{A\kon B,C,D}\cirk
b^\rightarrow_{A,B,C\kon D}=(b^\rightarrow_{A,B,C}\kon{\bf
1}_D)\cirk b^\rightarrow_{A,B\kon C,D}\cirk({\bf 1}_A\kon
b^\rightarrow_{B,C,D})$.
\end{tabbing}

\noindent We also assume besides reflexivity, symmetry and
transitivity of equality that if ${f=g}$ and ${h=j}$, then
${f\cirk h=g\cirk j}$, provided ${f\cirk h}$ and ${g\cirk j}$ are
defined, and ${f\kon h=g\kon j}$.

In \mak\ we have that $\kon$ is a bifunctor, $b^\rightarrow$ is a
natural isomorphism in all its indices, and \mbox{($b$5)} is Mac
Lane's \emph{pentagonal} equation of \cite{ML63}, where it is
proved that the category \mak\ is a preorder. Namely, for all
arrows ${f,g\!:\! A\rightarrow B}$ of \mak\ we have that ${f=g}$
(see also \cite{ML71}, Section VII.2, or \cite{DP04}, Section
4.3). The category \mak\ is the free monoidal category without
unit, i.e.\ free \emph{associative} category in the terminology of
\cite{DP04} (Section 4.3), generated by a single object, this
object being conceived as a trivial discrete category.

\section{The isomorphism of \G\ and \mak}

We are going to prove that the categories \G\ and \mak\ are
isomorphic. We define first what is missing of the structure of
\mak\ in \G\ in the following manner:

\begin{tabbing}
\mbox{\hspace{9em}}\= $b^\rightarrow_{A,B,C}$ \= $=_{\df}$ \=
$((\g22\triangleleft_3{\bf 1}_C)\triangleleft_2{\bf
1}_B)\triangleleft_1{\bf 1}_A$,
\\[2ex]
\> $b^\leftarrow_{A,B,C}$ \> $=_{\df}$ \>
$((\gamma^\leftarrow_{{\bf 2},{\bf 2}}\triangleleft_3{\bf
1}_C)\triangleleft_2{\bf 1}_B)\triangleleft_1{\bf 1}_A$,
\\[2ex]
\> \> $f\kon g$\' $=_{\df}$ \> $({\bf 1}_{\bf 2}\triangleleft_2
g)\triangleleft_1 f$.
\end{tabbing}

\noindent It can then be checked by induction on the length of
derivation that the equations of \mak\ are satisfied in \G.

We have of course the equations \mbox{(\emph{cat}~1)} and
\mbox{(\emph{cat}~2)}, while the equations \mbox{(\emph{bif}~1)}
and \mbox{(\emph{bif}~2)} with $\triangleleft_n$ replaced by
$\kon$ are easy consequences of \mbox{(\emph{bif}~1)} and
\mbox{(\emph{bif}~2)}. To derive \mbox{($b$~\emph{nat})}, we have
that with the help of \mbox{(\emph{assoc}~1$\rightarrow$)} and
\mbox{(\emph{bif}~1)} the left-hand side is equal to

\[
(((\g22\triangleleft_3{\bf 1}_F)\triangleleft_2{\bf
1}_D)\triangleleft_1{\bf 1}_B)\cirk((({\bf 1}_{{\bf
2}\triangleleft_2{\bf 2}}\triangleleft_3 h)\triangleleft_2
g)\triangleleft_1 f),
\]

\noindent while with the help of
\mbox{(\emph{assoc}~1$\rightarrow$)},
\mbox{(\emph{assoc}~2$\rightarrow$)} and \mbox{(\emph{bif}~1)} the
right-hand side is equal to

\[
((({\bf 1}_{{\bf 2}\triangleleft_1{\bf 2}}\triangleleft_3
h)\triangleleft_2 g)\triangleleft_1
f)\cirk(((\g22\triangleleft_3{\bf 1}_E)\triangleleft_2{\bf
1}_C)\triangleleft_1{\bf 1}_A).
\]

\noindent Then it is enough to apply \mbox{(\emph{bif}~2)} and
\mbox{(\emph{cat}~1)}. It is trivial to derive \mbox{($bb$)} with
the help of \mbox{(\emph{bif}~2)}, \mbox{($\gamma\gamma$)} and
\mbox{(\emph{bif}~1)}.

We derive finally the pentagonal equation \mbox{($b$5)}. With the
help of \mbox{(\emph{bif}~1)},
\mbox{(\emph{assoc}~1$\rightarrow$)} and
\mbox{(\emph{assoc}~2$\rightarrow$)} we derive that each of

\[
b^\rightarrow_{A,B,C\kon D},\;\;\; b^\rightarrow_{A\kon
B,C,D},\;\;\; {\bf 1}_A\kon b^\rightarrow_{B,C,D},\;\;\;
b^\rightarrow_{A,B\kon C,D},\;\;\; b^\rightarrow_{A,B,C}\kon{\bf
1}_D
\]

\noindent is equal to $(((f\triangleleft_4{\bf
1}_D)\triangleleft_3{\bf 1}_C)\triangleleft_2{\bf
1}_B)\triangleleft_1{\bf 1}_A$ for $f$ being respectively

\[
\g22\triangleleft_3{\bf 1}_{\bf 2},\;\;\; \g22\triangleleft_1{\bf
1}_{\bf 2},\;\;\; {\bf 1}_{\bf 2}\triangleleft_2\g22,\;\;\;
\g22\triangleleft_2{\bf 1}_{\bf 2},\;\;\; {\bf 1}_{\bf
2}\triangleleft_1\g22.
\]

\noindent Then, by relying on \mbox{(\emph{bif}~2)}, it is enough
to derive the following:

\begin{tabbing}

\mbox{\hspace{1em}}$(\g22\triangleleft_1{\bf 1}_{\bf
2})\cirk(\g22\triangleleft_3{\bf 1}_{\bf 2}) =
\gamma^\rightarrow_{{\bf 2}\triangleleft_1{\bf 2},{\bf
2}}\cirk(\g22\triangleleft_3{\bf 1}_{\bf 2})$, \quad by
\mbox{(\emph{hex}~1\emph{a})},
\\*[1ex]
\mbox{\hspace{9em}}\= = \= $({\bf 1}_{\bf
2}\triangleleft_1\g22)\cirk\gamma^\rightarrow_{{\bf
2}\triangleleft_2{\bf 2},{\bf 2}}$, \quad by
\mbox{($\gamma$~\emph{nat})},
\\[1ex]
\> = \> $({\bf 1}_{\bf
2}\triangleleft_1\g22)\cirk(\g22\triangleleft_2{\bf 1}_{\bf
2})\cirk({\bf 1}_{\bf 2}\triangleleft_2\g22)$, \quad by
\mbox{(\emph{hex}~1)}.

\end{tabbing}

\noindent Diagrammatically, we have

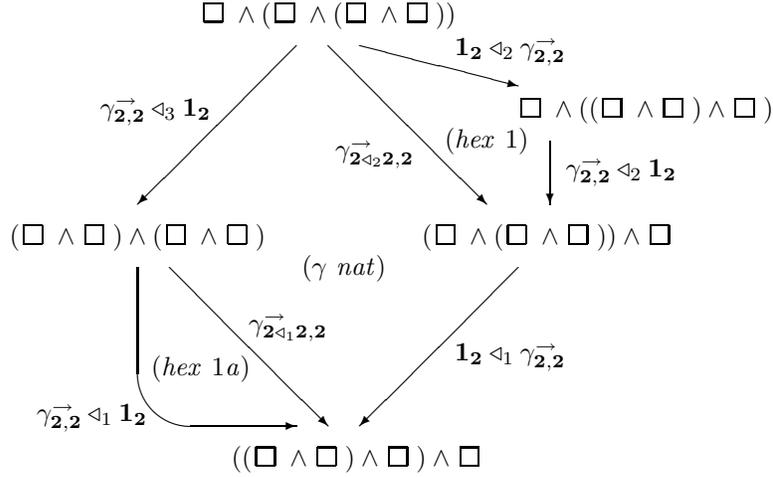
\begin{center}
\begin{picture}(180,175)
\unitlength1.2pt

\put(10,60){\vector(1,-1){50}} \put(120,60){\vector(-1,-1){50}}
\put(50,130){\vector(-1,-1){50}} \put(60,130){\vector(1,-1){50}}
\put(130,100){\vector(0,-1){20}} \put(70,130){\vector(4,-1){50}}
\put(40,10){\vector(1,0){10}}

\put(50,60){\oval(100,100)[bl]}

\put(70,0){\makebox(0,0){$((\koc\kon\koc)\kon\koc)\kon\koc$}}
\put(60,140){\makebox(0,0){$\koc\kon(\koc\kon(\koc\kon\koc))$}}
\put(0,70){\makebox(0,0){$(\koc\kon\koc)\kon(\koc\kon\koc)$}}
\put(130,70){\makebox(0,0){$(\koc\kon(\koc\kon\koc))\kon\koc$}}
\put(160,110){\makebox(0,0){$\koc\kon((\koc\kon\koc)\kon\koc)$}}

\put(20,28){\makebox(0,0){\mbox{(\emph{hex}~1\emph{a})}}}
\put(110,100){\makebox(0,0){\mbox{(\emph{hex}~1)}}}
\put(65,60){\makebox(0,0){\mbox{($\gamma$~\emph{nat})}}}
\put(100,127){\makebox(0,0)[bl]{${\bf 1}_{\bf
2}\triangleleft_2\g22$}}

\put(135,90){\makebox(0,0)[l]{$\g22\triangleleft_2{\bf 1}_{\bf
2}$}}

\put(87,100){\makebox(0,0)[tr]{$\gamma^\rightarrow_{{\bf
2}\triangleleft_2{\bf 2},{\bf 2}}$}}

\put(100,36){\makebox(0,0)[tl]{${\bf 1}_{\bf
2}\triangleleft_1\g22$}}

\put(35,40){\makebox(0,0)[bl]{$\gamma^\rightarrow_{{\bf
2}\triangleleft_1{\bf 2},{\bf 2}}$}}

\put(3,17){\makebox(0,0)[tr]{$\g22\triangleleft_1{\bf 1}_{\bf
2}$}}

\put(23,108){\makebox(0,0)[br]{$\g22\triangleleft_3{\bf 1}_{\bf
2}$}}
\end{picture}
\end{center}

\vspace{2ex}

\noindent So the pentagon is decomposed into a triangle (a
degenerate hexagon, corresponding to \mbox{(\emph{hex}~1)}), a
square (analogous to a naturality square, corresponding to
\mbox{($\gamma$~\emph{nat})}) and a two-sided diagram
(corresponding to \mbox{(\emph{hex}~1\emph{a})}).

If ${c_{A,B}\!:A\kon B\rightarrow B\kon A}$ is the commutativity
arrow of symmetric monoidal categories, for which in strict
categories of this kind, where associativity arrows are
identities, we have the equations

\begin{tabbing}

\hspace{4em}\=($c$ \emph{nat})\hspace{3em}\=$c_{B,D}\cirk(f\kon
g)=(g\kon f)\cirk c_{A,B}$,
\\[1.5ex]
\>($c$~\emph{hex}~1)\>$c_{A\kon B,C}=(c_{A,C}\kon{\bf
1}_B)\cirk({\bf 1}_A\kon c_{B,C})$,
\end{tabbing}

\noindent then we derive the Yang-Baxter equation

\begin{tabbing}

\hspace{0em}\=$(c_{B,C}\kon\!{\bf 1}_A)\cirk({\bf 1}_B\!\kon
c_{A,C})\cirk(c_{A,B}\kon\!{\bf 1}_C)=({\bf 1}_C\kon\!
c_{A,B})\cirk(c_{A,C}\kon\!{\bf 1}_B)\cirk({\bf 1}_A\!\kon
c_{B,C})$
\end{tabbing}

\noindent in the following way:

\begin{tabbing}
\hspace{1em}\=$(c_{B,C}\kon\!{\bf 1}_A)\cirk({\bf 1}_B\!\kon
c_{A,C})\cirk(c_{A,B}\kon\!{\bf 1}_C)=c_{B\kon
A,C}\cirk(c_{A,B}\kon\!{\bf 1}_C)$,\hspace{1em}by
($c$~\emph{hex}~1),
\\[1ex]
\hspace{10em}\=$=({\bf 1}_C\!\kon c_{A,B})\cirk c_{A\kon
B,C}$,\hspace{1em}by ($c$ \emph{nat}),
\\[1ex]
\>$=({\bf 1}_C\kon\! c_{A,B})\cirk(c_{A,C}\kon\!{\bf
1}_B)\cirk({\bf 1}_A\!\kon c_{B,C})$,\hspace{1em}by
($c$~\emph{hex}~1).
\end{tabbing}

\noindent This derivation is analogous to our derivation of
\mbox{($b$5)} above, where however the arrow corresponding to
${\bf 1}_B\!\kon c_{A,C}$ is identity, in virtue of the equation
\mbox{(\emph{assoc} 2)} on objects.

Alternatively, we derive \mbox{($b$5)} by using the following:

\begin{tabbing}

\mbox{\hspace{1em}}$(\g22\triangleleft_1{\bf 1}_{\bf
2})\cirk(\g22\triangleleft_3{\bf 1}_{\bf 2}) =
(\g22\triangleleft_1{\bf 1}_{\bf 2})\cirk\gamma^\rightarrow_{{\bf
2},{\bf 2}\triangleleft_2{\bf 2}}$, \quad by
\mbox{(\emph{hex}~2\emph{a})},
\\*[1ex]
\mbox{\hspace{9em}}\= = \= $\gamma^\rightarrow_{{\bf 2},{\bf
2}\triangleleft_1{\bf 2}}\cirk({\bf 1}_{\bf
2}\triangleleft_2\g22)$, \quad by \mbox{($\gamma$~\emph{nat})},
\\[1ex]
\> = \> $({\bf 1}_{\bf
2}\triangleleft_1\g22)\cirk(\g22\triangleleft_2{\bf 1}_{\bf
2})\cirk({\bf 1}_{\bf 2}\triangleleft_2\g22)$, \quad by
\mbox{(\emph{hex}~2)}.

\end{tabbing}

\noindent Diagrammatically, we have

\begin{center}
\begin{picture}(180,175)
\unitlength1.2pt

\put(0,60){\vector(1,-1){50}} \put(110,60){\vector(-1,-1){50}}
\put(60,130){\vector(-1,-1){50}} \put(70,130){\vector(1,-1){50}}
\put(120,20){\vector(-4,-1){50}} \put(130,60){\vector(0,-1){20}}
\put(0,90){\vector(0,-1){10}}

\put(50,80){\oval(100,100)[tl]}

\put(60,0){\makebox(0,0){$((\koc\kon\koc)\kon\koc)\kon\koc$}}
\put(60,140){\makebox(0,0){$\koc\kon(\koc\kon(\koc\kon\koc))$}}
\put(0,70){\makebox(0,0){$(\koc\kon\koc)\kon(\koc\kon\koc)$}}
\put(130,70){\makebox(0,0){$\koc\kon((\koc\kon\koc)\kon\koc)$}}
\put(160,30){\makebox(0,0){$(\koc\kon(\koc\kon\koc))\kon\koc$}}

\put(20,113){\makebox(0,0){\mbox{(\emph{hex}~2\emph{a})}}}
\put(105,35){\makebox(0,0){\mbox{(\emph{hex}~2)}}}
\put(65,80){\makebox(0,0){\mbox{($\gamma$~\emph{nat})}}}
\put(97,108){\makebox(0,0)[bl]{${\bf 1}_{\bf
2}\triangleleft_2\g22$}}

\put(135,50){\makebox(0,0)[l]{$\g22\triangleleft_2{\bf 1}_{\bf
2}$}}

\put(90,45){\makebox(0,0)[br]{$\gamma^\rightarrow_{{\bf 2},{\bf
2}\triangleleft_1{\bf 2}}$}}

\put(100,12){\makebox(0,0)[tl]{${\bf 1}_{\bf
2}\triangleleft_1\g22$}}

\put(37,102){\makebox(0,0)[tl]{$\gamma^\rightarrow_{{\bf 2},{\bf
2}\triangleleft_2{\bf 2}}$}}

\put(4,128){\makebox(0,0)[br]{$\g22\triangleleft_3{\bf 1}_{\bf
2}$}}

\put(23,33){\makebox(0,0)[tr]{$\g22\triangleleft_1{\bf 1}_{\bf
2}$}}
\end{picture}
\end{center}

\vspace{2ex}

\noindent This is an alternative decomposition of the pentagon
into a triangle, a square and a two-sided diagram. Hence we have
in \G\ all the equations of \mak.

To define what is missing of the structure of \G\ in \mak, we have
first the following inductive definition of $\triangleleft_n$ on
arrows:

\begin{tabbing}

\hspace{1em}if $A'\kon(B'\kon C')=(A\kon(B\kon C))\triangleleft_n
D$,
\\[1.5ex]
\hspace{3em}\=$b^\rightarrow_{A,B,C}\triangleleft_n{\bf
1}_D=b^\rightarrow_{A',B',C'}$,\quad\quad
$b^\leftarrow_{A,B,C}\triangleleft_n{\bf
1}_D=b^\leftarrow_{A',B',C'}$,
\\[2ex]
\>$(g\cirk f)\triangleleft_n{\bf 1}_D=(g\triangleleft_n{\bf
1}_D)\cirk(f\triangleleft_n{\bf 1}_D)$,
\\[2ex]
\>$(f\kon g)\triangleleft_n {\bf 1}_D=\left\{\begin{array}{ll}
(f\triangleleft_n {\bf 1}_D)\kon g & \mbox{\rm{if }} 1\leq n\leq
|f|
\\[1.5ex]
f\kon(g\triangleleft_{n-|f|} {\bf 1}_D) & \mbox{\rm{if }}
|f|<n\leq |f|+|g|,
\end{array}\right.$
\\[3ex]
\>${\bf 1}_{\koci}\triangleleft_1 f=f$,
\\*[2ex]
\>${\bf 1}_{A\kon B}\triangleleft_n f=\left\{\begin{array}{ll}
({\bf 1}_A\triangleleft_n f)\kon {\bf 1}_B &\;\; \mbox{\rm{if }}
1\leq n\leq |A|
\\[1.5ex]
{\bf 1}_A\kon({\bf 1}_B\triangleleft_{n-|A|} f) &\;\; \mbox{\rm{if
}} |A|<n\leq |A|+|B|,
\end{array}\right.$
\\[3ex]
\>$f\triangleleft_n g=({\bf 1}_B\triangleleft_n
g)\cirk(f\triangleleft_n{\bf 1}_C)$.

\end{tabbing}

We define $\gamma^\rightarrow_{A,B}$ and $\gamma^\leftarrow_{A,B}$
by stipulating

\[
\g22=_{\df}\;b^\rightarrow_{\koci,\koci,\koci},
\quad\quad\quad\gamma^\leftarrow_{{\bf 2},{\bf
2}}=_{\df}\;b^\leftarrow_{\koci,\koci,\koci},
\]

\noindent and by using \mbox{($\gamma${\bf 1})},
\mbox{(\emph{hex}~1)}, \mbox{(\emph{hex}~1\emph{a})},
\mbox{(\emph{hex}~2)} and \mbox{(\emph{hex}~2\emph{a})} as clauses
in an inductive definition.

The equations of \G\ certainly hold in \mak\ for this defined
structure because \mak\ is a preorder, as we said above. To finish
the proof that \G\ and \mak\ are isomorphic categories, it remains
only to check that the clauses of the inductive definitions of
$\triangleleft_n$, $\gamma^\rightarrow_{A,B}$ and
$\gamma^\leftarrow_{A,B}$ hold as equations in \G\ for
$b^\rightarrow_{A,B,C}$, $b^\leftarrow_{A,B,C}$ and $\kon$ defined
as they are defined in \G. This is done by using essentially
\mbox{(\emph{assoc}~1$\rightarrow$)} and
\mbox{(\emph{assoc}~2$\rightarrow$)}. So \G\ is isomorphic to
\mak, and is hence a preorder.

If we have instead of \mak\ the free monoidal category without
unit, i.e.\ the free associative category,
$\raisebox{2pt}{\makebox(8,5){\mak}}'$ generated by an arbitrary
nonempty set of objects $\cal P$, conceived as a discrete
category, then, instead of \G, the analogous category ${\bf
\Gamma}'$ isomorphic to $\raisebox{2pt}{\makebox(8,5){\mak}}'$
would have as generators ${\cal P}\cup\{{\bf 2}\}$. Every object
of ${\bf \Gamma}'$ different from a member of $\cal P$ can be
written in the form

\[
(\ldots(C\triangleleft_n p_n)\ldots\triangleleft_2
p_2)\triangleleft_1 p_1
\]

\noindent for $C$ an object of \G\ (more precisely, a member of
$\eL_2$), ${n=|C|}$ and ${p_1,\ldots,p_n}$ members of $\cal P$.
For the arrows $\gamma^\rightarrow_{A,B}$ and
$\gamma^\leftarrow_{A,B}$ we would assume that $p_{|A|}$ in $A$
coincides with $p_1$ in $B$, and the equations \mbox{($\gamma${\bf
1})} and \mbox{(\emph{unit} $\rightarrow$)} would have to be
adapted.

We have seen above how Mac Lane's pentagon arises from a
Yang-Baxter hexagon by collapsing, according to
\mbox{(\emph{assoc}~2)}, the vertices corresponding to ${({\bf
2}\triangleleft_1{\bf 2})\triangleleft_3{\bf 2}}$, i.e.\ ${B\kon
A\kon C}$, and ${({\bf 2}\triangleleft_2{\bf
2})\triangleleft_1{\bf 2}}$, i.e.\ ${B\kon C\kon A}$, into a
single vertex corresponding to
${(\koc\kon\koc)\kon(\koc\kon\koc)}$. We can apply this collapsing
procedure based on \mbox{(\emph{assoc}~2)} to the
three-dimensional permutohedron (whose vertices correspond to
permutations of four letters and edges to transpositions of
adjacent letters) in order to obtain the three-dimensional
associahedron (whose vertices correspond to planar binary trees
with five leaves and edges to arrow terms of \mak\ with a single
$b^\rightarrow$), and afterwards we can proceed to higher
dimensions. The function that corresponds to our procedure is
described in \cite{T97}. Our paper provides a motivation for that
function.

\baselineskip=0.84\baselineskip

\end{document}